\documentclass[12pt,reqno]{amsart}
\usepackage{amsmath}
\usepackage{amssymb}
\usepackage{amsthm}
\usepackage[left=2.9cm,top=2.5cm,right=2.9cm,bottom=2.6cm]{geometry}
\usepackage{epsfig}
\usepackage[cp1250]{inputenc}
\usepackage[T1]{fontenc}

\begin{document}
\newtheorem{theorem}{Theorem}[section]
\newtheorem{lemma}[theorem]{Lemma}
\newtheorem{definition}[theorem]{Definition}
\newtheorem{conjecture}[theorem]{Conjecture}
\newtheorem{proposition}[theorem]{Proposition}
\newtheorem{algorithm}[theorem]{Algorithm}
\newtheorem{corollary}[theorem]{Corollary}
\newtheorem{observation}[theorem]{Observation}
\newtheorem{problem}[theorem]{Problem}
\newtheorem{example}[theorem]{Example}
\newcommand{\noin}{\noindent}
\newcommand{\ind}{\indent}
\newcommand{\al}{\alpha}
\newcommand{\om}{\omega}
\newcommand{\pp}{\mathcal P}
\newcommand{\ppp}{\mathfrak P}
\newcommand{\R}{{\mathbb R}}
\newcommand{\N}{{\mathbb N}}
\newcommand{\Z}{{\mathbb Z}}
\newcommand\eps{\varepsilon}
\newcommand{\E}{\mathbb E}
\newcommand{\Prob}{\mathbb{P}}
\newcommand{\pl}{\textrm{C}}
\newcommand{\dang}{\textrm{dang}}
\renewcommand{\labelenumi}{(\roman{enumi})}
\newcommand{\bc}{\bar c}
\newcommand{\G}{{\mathfrak S}}
\newcommand{\T}{{\mathfrak T}}
\newcommand{\Remark}[1]{\par \noindent \textnormal{\textbf{Note.~}} #1}

\title{Planar graph is on fire}
\author{Przemys\l{}aw Gordinowicz}
\address{Institute of Mathematics, Technical University of Lodz, \L{}\'od\'z, Poland}
\email{pgordin@p.lodz.pl}
\begin{abstract}
Let $G$ be any connected graph on $n$ vertices, $n \ge 2.$ Let $k$ be any positive integer. Suppose that a fire breaks out at some vertex of $G.$ Then, in each turn $k$ firefighters can protect vertices of $G$ --- each can protect one vertex not yet on fire; Next the fire spreads to all unprotected neighbours. 

The \emph{$k$-surviving} rate of G, denoted by  $\rho_k(G),$ is the expected fraction of vertices that can be saved from the fire by $k$ firefighters, provided that the starting vertex is chosen uniformly at random. In this paper, it is shown that for any planar graph $G$ we have $\rho_3(G) \ge \frac{2}{21}.$ Moreover, 3 firefighters are needed for the first step only; after that it is enough to have 2 firefighters per each round. This result significantly improves the known solutions to a problem by Cai and Wang (there was no positive bound known for the surviving rate of general planar graph with only 3 firefighters). The proof is done using the separator theorem for planar graphs.
\end{abstract}
\keywords{firefighter problem, surviving rate, planar graph}
\subjclass[2000]{05C15}


\vspace{150pt}

\maketitle

\section{Introduction\label{sec:intro}}
The following \emph{Firefighter Problem} was introduced by Hartnell~\cite{hartnell}. Consider any connected graph, say $G,$ on $n$ vertices, $n \ge 2.$ Let $k$ be any positive integer. Suppose that a fire breaks out at some vertex $v \in V(G).$ Then in each turn firefighters can protect $k$ vertices of $G,$ not yet on fire and the protection is permanent. Next the fire spreads to all the unprotected vertices that are adjacent to some vertices already on fire. The goal is to save as much as possible and the question is how many vertices can be saved. We would like to refer the reader to the survey by Finbow and MacGillivray~\cite{fm} for more information on the background of the problem and directions of its consideration. 

In this paper we focus on the following aspect of the problem. 
Let $sn_k(G, v)$ denote the maximum number of vertices of $G$ that $k$ firefighters can save when the fire breaks out at the vertex $v.$ This parameter may depend heavily on the choice of the starting vertex $v,$ for example when the graph $G$ is a star. Therefore Cai and Wang~\cite{cai} introduced the following graph parameter: the \emph{$k$-surviving rate} $\rho_{k}(G)$ is the expected fraction of vertices that can be saved by $k$ firefighters, provided that the starting vertex is chosen uniformly at random. Namely $$\rho_k(G) = \frac{1}{|V(G)|^2} \sum_{v \in V(G)} sn_k(G, v).$$

While discussing the surviving rate, let us mention the recent results by Prałat~\cite{pralat3011,pralatfire} which have provided a~threshold for the average degree of general graphs, which guarantees a positive surviving rate with a given number of firefighters. To be more precise, for $k \in \N^+$ let us define
$$\tau_k = \left\{ \begin{array}{lcl}
\frac{30}{11} & \textrm{~for~} & k = 1\\
k + 2 - \frac{1}{k+2}& \textrm{~for~} & k \ge 2.\\    
\end{array}\right.$$
Then, there exists a constant $c > 0,$ such that for any $\epsilon > 0,$ any $n \in \N^+$ and any graph $G$ on $n$ vertices and at most $(\tau_k - \epsilon)n/2$ edges one has $\rho_k(G) > c \cdot \epsilon >0.$ Moreover, there exists a~family of graphs with the average degree tending to $\tau_k$ and the $k$-surviving rate tending to 0, which shows that the above result is the best possible.

The $k$-surviving rate is investigated for many particular families of graphs --- we focus here on planar graphs. Cai and Wang~\cite{cai} asked about the minimum number of firefighters $k$ such that $\rho_k(G) > c$ for some positive constant $c$ and any planar graph $G.$ It is easy to see that $\rho_1(K_{2,n}) \overrightarrow{_{\ n \to \infty}\ } \, 0,$ hence at least 2 firefighters are necessary. It is shown that 2 is the upper bound for triangle-free planar graphs~\cite{esperet} and planar graphs without 4-cycles~\cite{kwzhang}.

So far, for the general planar graphs the best known upper bound for the number of firefighters is 4: Kong, Wang and Zhu~\cite{kwzhu} have shown that $\rho_4(G) > \frac{1}{9}$ for any planar graph $G.$ Esperet, van den Heuvel, Maffray and Sipma~\cite{esperet} have shown that using 4 firefighters in the first round only and just 3 in the subsequent rounds it is also possible to save a positive fraction of any planar graph $G,$ namely $\rho_{4,3}(G) > \frac{1}{2712}.$ We use the notation of $\rho_{k,l}$ and $sn_{k,l}$ to describe the model with $k$ firefighters in the first round and $l$ firefighters in the subsequent rounds.

In this paper, we improve the above bounds by the following theorem:
\begin{theorem} \label{thm:main}
Let $G$ be any planar graph. Then:
\begin{enumerate}
\item $\rho_{4,2}(G) > \frac{2}{9}.$
\item $\rho_{3}(G) \ge \rho_{3,2}(G) > \frac{2}{21}.$
\end{enumerate}
\end{theorem}
In other words, we show that with 3 firefighters in the first round and just 2 in the subsequent rounds we can save at least $\frac{2}{21}$ vertices of a planar graph, while with one extra firefighter in the first round we can increase the saved fraction to $\frac{2}{9}.$ 

\section{The proof\label{sec:proof}}
The proof is done using the lemma given by Lipton and Tarjan to prove the separator theorem for planar graphs~\cite{liptontarjan}. The key lemma in their proof, slightly reformulated to use in the firefighter problem, is quoted below.
\begin{lemma} \label{thm:separator}
Let $G$ be any $n$-vertex plane triangulation and $T$ be any spanning tree of $G.$ Then there exists an edge $e \in E(G) \setminus E(T)$ such that the only cycle $C$ in $T + e$ has the property that the number of vertices inside $C$ as well as outside $C$ is lower than $\frac{2}{3}n.$
\end{lemma}
A similar approach --- using the above lemma to the firefighter problem on planar graphs, was first applied by Floderus, Lingas and Persson~\cite{lingas2}, with a slightly different notation of approximation algorithms. The authors of~\cite{lingas2} have proved a theorem analogous to Lemma~\ref{lem:many}. For some more details see Section~\ref{sec:remarks}.

The proof of Theorem~\ref{thm:main} is presented in two steps --- first we show that $\rho_{4,2}(G) > \frac{2}{9}$ for any planar graph $G$, then that $\rho_{3,2}(G) > \frac{2}{21}.$ At first let us note that the surviving rate is monotone (non-increasing) with respect to the operation of adding edges to the graph. Hence, it is enough to prove the bounds given by Theorem~\ref{thm:main} only for plane triangulations. Moreover, in the first step, depending on the number of firefighters, we save 3 or 4 vertices respectively, which is enough to obtain the desired bounds for any planar graph on not more than 17 vertices. 

Let $G$ be any $n$-vertex plane triangulation, where $n \ge 18.$ Suppose that the fire breaks out at a vertex $r.$ Consider a tree $T$ obtained by the breadth-first-search algorithm starting from the vertex $r.$ By Lemma~\ref{thm:separator} there is the edge $e$ and the cycle $C \subseteq T + e$ such that $|C \cup \textrm{in} C| > \frac{1}{3}n$ and $|C \cup \textrm{out} C| > \frac{1}{3}n,$ where $\textrm{in} C$ and $\textrm{out} C$ denote the sets of vertices inside the cycle $C$ and outside the cycle $C$ respectively. Note that in the cycle $C$ there are at most 2 vertices at any given distance from $r.$ This holds because every edge of $C$, except one, belongs to the breadth-first-search tree. The firefighters' strategy depends on the cycle $C.$ When the vertex $r$ does not belong to the cycle then the firefighters protect the vertices of $C$ in the order given by the distance from the vertex $r$ and save all the vertices in either $C \cup \textrm{in} C$ or $C \cup \textrm{out} C.$ When the vertex $r$ belongs to the cycle $C,$ the firefighters still can protect the vertices of the cycle except the vertex $r,$ but it may be not enough, as the fire may spread through the neighbours of $r$ inside as well as outside the cycle $C.$
Because either $\textrm{in} C$ or $\textrm{out} C$ contains not more than $\left\lfloor\frac{\deg r -2}{2}\right\rfloor$ neighbours of $r$ we get immediately:
\begin{lemma} \label{lem:many}
Let $G$ be any $n$-vertex plane triangulation, where $n \ge 18.$ Suppose that the fire breaks out at some vertex $r.$ Then using $2 + \left\lfloor\frac{\deg r -2}{2}\right\rfloor$ firefighters at the first step and 2 at the subsequent steps one can save more than $n/3 - 1$ vertices. 
\end{lemma}
To calculate the surviving rate $\rho_{4,2}(G)$ let us now partition the vertex set of the graph $G$ into 3 sets: $X = \{v \in V(G) \colon \deg(v) \in \{3,4\}\},$ $Y =  \{v \in V(G) \colon \deg(v) \in \{5,6,7\}\}$ and $Z = \{v \in V(G) \colon \deg(v) \ge 8\}.$ Obviously we have $|Z| = n - |X| - |Y|.$ Since for the plane triangulation we have $$6n > 6n - 12 = \sum_{v \in V(G)} \deg v \ge 3|X|+5|Y|+8|Z|,$$ then one has $$|Y| > \frac{2n-5|X|}{3}.$$
Using 4 firefighters at the first step one can save $n-1$ vertices if $r \in X,$ more than $n/3 - 1$ vertices if $r \in Y$ and at least $4$ vertices if $r \in Z.$  A simple calculation shows now that $$\rho_{4,2}(G) > \frac{(n-1)|X|+(n/3-1)|Y|+4|Z|}{n^2} > \frac{2}{9},$$  
which finishes the proof for the case with 4 firefighters.

Let us start our proof for 3 firefighters with a simple observation derived from Lemma~\ref{lem:many}.
\begin{observation}
Let $G$ be any $n$-vertex plane triangulation, where $n \ge 18.$ Let $r \in V(G)$ be a vertex of degree at most $5.$ Then 
$$sn_{3,2}(G, r) > \left\{ \begin{array}{lcl}
n-1 & \textrm{~if~} & \deg(r) \le 3\\
n/3-1 & \textrm{~if~} & \deg(r) \in \{4,5\}\\ 
\end{array}\right.$$
\end{observation}

Dealing with vertices of degree higher than 5 is a bit more complicated. Of course, the firefighters still can save at least 3 vertices, but frequently enough it is possible to save more.
\begin{lemma} \label{lem:3f}
Let $G$ be any $n$-vertex plane triangulation, where $n \ge 18.$ Let $r \in V(G)$ be a vertex of degree $6$ or $7.$ Then either 
$$sn_{3,2}(G, r) > \left\{ \begin{array}{lcl}
n/4-1 & \textrm{~if~} & \deg(r) = 6\\
n/6-1 & \textrm{~if~} & \deg(r) = 7\\    
\end{array}\right.$$
or the vertex $r$ has at least $2$ adjacent neighbours, say $u$ and $v,$ such that $sn_{3,2}(G, u) > n/3$ and $sn_{3,2}(G, v) > n/3.$
\end{lemma}

\begin{proof}
Let $G$ be any $n$-vertex plane triangulation, where $n \ge 18.$ Let $r \in V(G)$ be a vertex of degree 6 or 7. Consider a tree $T$ obtained by the breadth-first-search algorithm starting from the vertex $r.$ By Lemma~\ref{thm:separator} there is the edge $e$ and the cycle $C \subseteq T + e$ such that $|C \cup \textrm{in} C| > \frac{1}{3}n$ and $|C \cup \textrm{out} C| > \frac{1}{3}n.$ If the vertex $r$ has no more than one neighbour either inside the cycle $C$ or outside the cycle, then $sn_{3,2}(G, r) \ge n/3 -1.$ 

Without loss of generality we may assume then that the vertex $r$ has exactly 2 neighbours inside the cycle $C$ and at least 2 neighbours outside. When the vertex $r$ has degree 6 we may assume additionally, without loss of generality, that the number of vertices inside $C$ is not lower than the number of vertices outside, that is $|C \cup \textrm{in} C| > \frac{1}{2}n.$ Note that the terms ,,inside'' and ,,outside'' the cycle depend on a particular drawing of the triangulation.  

Let $u$ and $v$ be the neighbours of the vertex $r$ inside the cycle. Then one of the two following cases occur:
\begin{itemize}
\item[Case 1.] In the graph $G$ there exists a path from $u$ or $v$ to a vertex on the cycle containing vertices in increasing distance from $r.$
\item[Case 2.] There is no such path. 
\end{itemize}
\begin{figure}[htbp]
\begin{minipage}[b]{0.5\textwidth}
\begin{flushright}
{\epsfig{clip, viewport=24 11 359 395, figure=Separator, width=2in}}\\
\end{flushright}
\end{minipage}
\begin{minipage}[t]{0.38\textwidth}
\vspace{-50pt}
{\small
Solid edges are the edges of the spanning tree, bold if they belong to the cycle. Dashed edges are the edges forming the path (possibly but not necessarily belonging the tree).}
\end{minipage}
\caption{Ilustration of Case 1.}\label{fig1}
\end{figure}
The path described in the first case divides the cycle $C$ into two cycles $C'$ and $C''$ (see Figure 1), both of which have the properties that there are at most 2 vertices at any given distance from $r$ and there is at most one neighbour of the vertex $r$ inside the cycle. So, with 3 firefighters in the first round and 2 in the subsequent rounds one is able to save every vertex except $r$ from $C' \cup \textrm{in} C'$ or from $C'' \cup \textrm{in} C''.$ Choosing the larger piece it is possible to save at least half of the vertices from the set $C \cup \textrm{in} C$ (note that the vertices on the path count for both pieces). 

Considering the second case, note that as $G$ is a triangulation then $u$ and $v$ are adjacent. Suppose now that the fire breaks out at the vertex $u$ or $v$ instead of the vertex $r.$ The firefighters can now save all the vertices in $C \cup \textrm{out} C$ by protecting in the first round vertex $r$ and its neighbours on the cycle $C,$ and the vertices of $C$ ordered in increasing distance from $r$ in the subsequent rounds. The fire cannot reach the vertices on the cycle $C$ earlier than the firefighters protect them --- otherwise, there would exist a vertex on the cycle $C$ which is closer to $u$ or $v$ than to $r$ --- this guarantees the existence of a path described in Case 1.
\end{proof}

Let us now partition the vertex set of $G$ into 4 subsets defined by the conditions:
\begin{eqnarray}
\nonumber X &=& \{v \in V(G) \colon sn_{3,2}(G, v) > \frac{n}{3}-1\},\\
\nonumber Y &=& \{v \in V(G) \colon \deg(v) \le 7 \wedge sn_{3,2}(G, v) \le \frac{2n}{21}\},\\
\nonumber Z &=& \{v \in V(G) \colon \deg(v) \ge 8 \wedge sn_{3,2}(G, v) \le \frac{2n}{21}\},\\ 
\nonumber W &=& V(G) \setminus (X \cup Y \cup Z).
\end{eqnarray}
We have that 
\begin{equation}\label{eq:sr3}
\rho_{3,2}(G) \ge \frac{1}{n^2} \Big(|X| (n/3-1) + 3|Y|+3|Z| + |W|2n/21\Big).
\end{equation}
Note that $X$ contains every vertex of $G$ with degree lower than 6. As the average degree of vertex in a plane triangulation is lower than 6 then the set $X$ is nonempty and the average degree of vertex in $X$ is also lower than $6.$ 
Every vertex $v \in W$ has $\deg(v) \ge 6$ and $sn_{3,2}(G, v) > \frac{2n}{21}.$ By Lemma~\ref{lem:3f} every vertex in the set $Y$ has at least 2 adjacent neighbours in the set $X,$ hence we have 
$$|Y| \le \frac{\sum_{x \in X} (\deg x - 1)}{2} = \frac{\sum_{x \in X} \deg x}{2} - \frac{|X|}{2}.$$
As any vertex in the set $Z$ has degree at least 8, while the average degree in $G$ is lower than 6, we have
$$|Z| < \frac{\sum_{x \in X}(6-\deg(x))}{2}.$$ These inequalities, when added, yield $|Y|+|Z| < \frac{5}{2} |X|.$ Hence, if only $3|Y|+3|Z| \ge |X|,$  Inequality~(\ref{eq:sr3}) yields
$$\rho_{3,2}(G) > \frac{\frac{|X|}{3}+\frac{2|W|}{21}}{\frac{7}{2}|X|+|W|} = \frac{2}{21}.$$ Note that it may occur that $3|Y|+3|Z| < |X|$ --- this is an easier case as then
$$\rho_{3,2}(G) > \frac{\frac{n|X|}{3} - |X| +\frac{2n|W|}{21}}{n(\frac{4}{3}|X|+|W|)} > \frac{\frac{8|X|}{63}+\frac{2|W|}{21}}{\frac{4}{3}|X|+|W|} = \frac{2}{21}.$$

\section{Remarks}~\label{sec:remarks}
While investigating whether or not the separator theorem had been used for the firefighter problem the author found~\cite{lingas} (recently it was published in the journal version as~\cite{lingas2}). There a theorem analogous to Lemma~\ref{lem:many} is proved. Moreover,~\cite{lingas,lingas2} present a theorem which in our notation would give that for any planar graph $G$ and any vertex $r$ one has $sn_{2}(G, r) \ge \frac{n}{3\deg r}.$ Unfortunately, there is a serious error in the proof: it is implicitly assumed that adding some edge joining two vertices of some induced subgraph of the planar graph, which preserves planarity of the subgraph, should also preserve planarity of the whole graph (or, in other words, that the separator of the subgraph is also a separator of the whole graph). In our opinion, such a~result cannot be proved just by a simple application of the separator theorem. Hence, the problem to determine whether the $2$-surviving rate of the planar graph may be bounded by some positive constant remains still open.


{\footnotesize
\begin{thebibliography}{99}
\bibitem{cai} L.~Cai, W.~Wang, The surviving rate of a graph for the firefighter problem, \emph{SIAM J. Discrete Math.} \textbf{23} (2009) 1814-–1826.

\bibitem{esperet} L.~Esperet, J.~van~den~Heuvel, F.~Maffray, F.~Sipma, Fire containment in planar graphs \emph{J. Graph Theory} \textbf{73} (2013) 267--279.

\bibitem{fm} S.~Finbow, G.~MacGillivray, The firefighter problem: a survey of results, directions and questions, \emph{Australasian Journal of Combinatorics} \textbf{43} (2009) 57–-77.

\bibitem{lingas} P.~Floderus, A.~Lingas, M.~Persson, Towards more efficient infection and fire fighting, Proceedings of the Seventeenth Computing: The Australasian Theory Symposium (CATS 2011), Perth, Australia, 69--73.

\bibitem{lingas2} P.~Floderus, A.~Lingas, M.~Persson, Towards more efficient infection and fire fighting, International Journal of Foundations of Computer Science \textbf{24} (2013) 3-–14.

\bibitem{hartnell}  B.~Hartnell, Firefighter! An application of domination, Presentation at the 25th Manitoba Conference on Combinatorial Mathematics and Computing, University of Manitoba, Winnipeg, Canada, 1995.

\bibitem{kwzhu} J.~Kong, W.~Wang, X.~Zhu, The surviving rate of planar graphs, \emph{Theoret. Comput. Sci.}, \textbf{416} (2012) 65-–70.

\bibitem{kwzhang} J.~Kong, W.~Wang, L.~Zhang, The 2-surviving rate of planar graphs without 4-cycles, \emph{Theoret. Comput. Sci.} \textbf{457} (2012) 158--165.

\bibitem{liptontarjan} R.~J.~Lipton, R.~E.~Tarjan, A Separator Theorem
for Planar Graphs, \emph{SIAM Journal on Applied Mathematics} \textbf{36} (1979), 177-–189.

\bibitem{pralat3011} P.~Pralat, Graphs with average degree smaller than 30/11 burn slowly, \emph{Graphs and Combinatorics}, \textbf{30} (2014), 455--470. 

\bibitem{pralatfire} P.~Pralat, Sparse graphs are not flammable, \emph{SIAM Journal on Discrete Mathematics} \textbf{27} (2013), 2157--2166. 

\end{thebibliography}}
\end{document}